\documentclass[3p,article]{elsarticle}

\usepackage{lineno,hyperref}
\usepackage{physics}
\usepackage{bm} % Maakt \boldsymbol van amsmath callbaar met \bm 
\usepackage{flafter} % Figure only appears after referred to
\usepackage[pdftex,dvipsnames]{xcolor}  % Coloured text etc.
\usepackage{tikz}
\usepackage{amssymb}
\usepackage{subcaption}
\usepackage{caption}
\usepackage{hyperref}
\hypersetup{
    colorlinks=true,
    linkcolor=blue,
    urlcolor=red,
    linktoc=all,
    pdfauthor=Lotz et. al.
}
\usepackage{pgfplots}
\pgfplotsset{compat=1.17}
\usepgfplotslibrary{colorbrewer}
\usepgfplotslibrary{patchplots}
\usepackage{graphicx}
\usepackage{comment}
\usepackage{nicefrac}
\usepackage{pdflscape}
\usepackage{geometry}

\makeatletter
\def\ps@pprintTitle{%
  \let\@oddhead\@empty
  \let\@evenhead\@empty
  \def\@oddfoot{\reset@font\hfil\thepage\hfil}
  \let\@evenfoot\@oddfoot
}
\makeatother

\usepackage[numbers]{natbib}
\usepackage{lineno}[pagewise]

\usepackage{xargs}                      % Use more than one optional parameter in a new commands
\usepackage[colorinlistoftodos,prependcaption,textsize=footnotesize]{todonotes}
\newcommandx{\unsure}[2][1=]{\todo[linecolor=red,backgroundcolor=red!25,bordercolor=red,#1]{#2}}
\newcommandx{\change}[2][1=]{\todo[linecolor=blue,backgroundcolor=blue!25,bordercolor=blue,#1]{#2}}
\newcommandx{\info}[2][1=]{\todo[linecolor=green,backgroundcolor=green!25,bordercolor=green,#1]{#2}}
\newcommandx{\improve}[2][1=]{\todo[linecolor=orange,backgroundcolor=orange!25,bordercolor=orange,#1]{#2}}
\newcommandx{\thiswillnotshow}[2][1=]{\todo[disable,#1]{#2}}

\usepackage[noabbrev,capitalise,nameinlink]{cleveref} % clever references --> \cref --> Also includes "figure", also multiple options are possible \cref{a,b} \cpageref --> revers to page, noabbrev stops the abbreviation, nameinlink also make the "figure" clickable. SHOULD BE LOADED LAST
\crefformat{equation}{(#2#1#3)}

\modulolinenumbers[5]

\journal{Journal of \LaTeX\ Templates}

\newcommand{\bx}{\bm{\mathrm{x}}} % Bold x (space coordinate)
\newcommand{\bX}{\bm{\mathrm{X}}} % Bold X (space coordinate in Lagrangian frame of reference)
\newcommand{\bhatx}{\bm{\mathrm{\hat{x}}}} % Bold hat x (space-time coordinate)
\newcommand{\bxi}{\bm{\mathrm{\xi}}} % Bold xi (reference coordinate)
\newcommand{\st}{x_{d+1}} % Bold x (reference coordinate)

\newcommand{\bu}{\bm{\mathrm{u}}} % Bold u (space velocity)
\newcommand{\un}{u_\mathrm{n}} % Bold u (space velocity)
\newcommand{\sbu}{\bm{\mathrm{u'}}} % Bold u (small scales of velocity)
\newcommand{\bhatu}{\bm{\mathrm{\hat{u}}}} % Bold hat u (space-time velocity)
\newcommand{\bg}{\bm{\mathrm{g}}} % Bold g (boundary velocity)
\newcommand{\bO}{\bm{\mathrm{0}}} % Bold 0 (zero velocity vector)
\newcommand{\vst}{{s} } % Velocity relating space and time
\newcommand{\p}{p} % Scalar p (pressure)
\newcommand{\bfo}{\bm{\mathrm{f}}} % Bold f (force)

\newcommand{\bU}{\bm{{U}}} % Bold U variable space
\newcommand{\bW}{\bm{{W}}} % Bold W testspace

\newcommand{\nb}{\Gamma^N_{\textrm{ext}}} % Neumann boundary domain
\newcommand{\db}{\Gamma^D_{\textrm{ext}}} % Dirichlet boundary domain
\newcommand{\nsb}{\Gamma_\textrm{int}} % No slip boundary domain
\newcommand{\eb}{\Gamma_{\textrm{ext}}} % Exterior boundary
\newcommand{\nsbs}{P_\textrm{int}} % Space time no slip boundary domain

\newcommand{\nbs}{P^N_{\textrm{ext}}} % Space time Neumann boundary domain
\newcommand{\dbs}{P^D_{\textrm{ext}}} % Space time Dirichlet boundary domain

\newcommand{\bn}{\bm{\mathrm{n}}} % Normal
\newcommand{\bhatn}{\bm{\mathrm{\hat{n}}}} % Space-time normal

\newcommand{\bG}{\bm{\mathrm{G}}} % Bold G (Gij)
\newcommand{\bhatG}{\bm{\mathrm{\hat{G}}}} % Bold hat G (Gij)
\newcommand{\br}{\bm{\mathrm{r}}} % Residuals

\newcommand{\bw}{\bm{\mathrm{w}}} % Weight function
\newcommand{\gbw}{\nabla\bm{\mathrm{w}}} % Gradient of weight function
\newcommand{\divbw}{\nabla\cdot\bm{\mathrm{w}}} % Gradient of weight function
\newcommand{\gbwhat}{\nabla_\mathrm{{\bhatx}}\bm{\mathrm{w}}} % Gradient of weight function to hat x
\newcommand{\gbq}{\nabla q} % Gradient of weight function of mass cons

\newcommand{\WW}{\mathcal{W}}

\newcommand{\UU}{\mathcal{U}}
\newcommand{\PP}{\mathcal{P}}

\begin{document}
\begin{frontmatter}

\title{A space-time framework for periodic flows with applications to hydrofoils}

\author[1]{Jacob E. Lotz \corref{mycorrespondingauthor}}
\address[1]{Delft University of Technology, Department of Mechanical, Maritime and Materials Engineering. P.O. Box 5, 2600 AA Delft, The Netherlands}
\cortext[mycorrespondingauthor]{Corresponding author}
\ead{j.e.lotz@tudelft.nl}

\author[2]{Marco F.P. ten Eikelder}
\address[2]{Institute for Mechanics, Computational Mechanics Group, Technical University of Darmstadt, Franziska-Braun-Straße 7, 64287 Darmstadt, Germany}
\ead{marco.eikelder@tu-darmstadt.de}
\author[1]{Ido Akkerman}
\ead{i.akkerman@tudelft.nl}

\journal{Computers and Fluids}

\begin{abstract}
In this paper we propose a space-time framework for the computation of periodic flows. We employ the isogeometric analysis framework to achieve higher-order smoothness in both space and time. The discretization is performed using residual-based variational multiscale modelling and weak boundary conditions are adopted to enhance the accuracy near the moving boundaries of the computational domain. We show conservation properties and present a conservative method for force extraction. We apply our framework to the computation of a heaving and pitching hydrofoil. Numerical results display very accurate results on course meshes.
\end{abstract}


\end{frontmatter}

\section{Introduction}

Periodic flows are ubiquitous in a large number of industrial applications and natural features. Prototypical examples include the flow around submerged propellers, wind turbines, or rotating flows in turbomachines and engines and the pulsatile flow of blood. Various challenges arise in the design of practical numerical simulations of these flows. On top of the well-known complications centered around inertia-driven character and the imposition of boundary conditions, the periodic nature adds novel peculiar hurdles. The typical strategy of simulating a periodic flow problem is to perform an unsteady computation in which the flow slowly develops periodic characteristics \cite{Kinsey2008, Michelassi2003}. As such, the computational time far exceeds that of one single period. Moreover, a user-defined criterion of the characteristics of the flow is inevitable and the flow is never strictly periodic. In this work we exploit the periodic nature of the problem and propose a space-time finite element method in the framework of residual-based variational multiscale (VMS) methods, isogeometric analysis and weak boundary conditions. Particular emphasis is on the application to heaving and pitching hydrofoils.

The concept of \textit{space-time finite elements} may be traced back to the late sixties, with important contributions by Fried \cite{fried1969finite} and Oden \cite{oden1969general,oden1969general2}. In time-dependent problems, the standard is to separate the discretization of the time (e.g. finite difference schemes) and space (e.g. Galerkin methods). This is often referred to as the \textit{semi-discrete} method. The idea of space-time finite element methods is to adopt the variational approach in the space-time setting \cite{Hulbert1990, Mittal1992, Masud1997}. An overview of applications of the space-time method can be found in \cite{Tezduyar2019}. In the context of fluid mechanics, seminal contributions are the works by Hughes and collaborators on stabilized methods in the space-time framework \cite{hughes1989new,shakib1991new}. A few years later, the \textit{VMS framework} \cite{hughes1995multiscale,hughes1998variational,Yan2017}, encompassing many existing stabilized methods, was proposed. The framework was originally introduced for stationary problems. In \cite{hughes1996space} it was argued that the most theoretically coherent framework for the extension to time-dependent problems is the space-time context. The most popular applications of the VMS methodology for time-dependent problems are however in the semi-discrete setting. A notable contribution in this regard is the work \cite{Bazilevs2007} that presented a variationally consistent VMS methodology for turbulent flows called \textit{residual-based variational multiscale} (RBVMS). This method is often used in combination with weak boundary conditions \cite{Bazilevs2007a}. 
Recently the popular Nitsche's method for the imposition of weak boundary conditions has been identified as a variational multiscale formulation \cite{stoter2021nitsche}. 
RBVMS opened the door for the development of a novel class of small-scale models for large-eddy simulations, including dynamic small-scales \cite{codina2002stabilized,EiAk17ii,evans2020variational} and discontinuity capturing \cite{ten2020theoretical,ten2019variation}.
The last important development with implications for the space-time framework that we succinctly discuss is the introduction of isogeometric analysis \cite{Hughes2005,cottrell2009isogeometric}. In contrast to classical space-finite element methods, isogeometric analysis offers the possibility of arbitrary smooth finite element basis functions. This technique was initially adopted for spatial discretizations, yet it offers rich opportunities in the space-time setting \cite{Takizawa2014,Kuraishi2019,Otoguro2017,montardini2020space,Saade2021}. On top of the well-known advantages of isogeometric analysis, the adoption of it in the space-time context is particularly beneficial for an accurate representation of moving boundaries and a higher continuity in the temporal direction.

The existing space-time finite element methods form a versatile and fundamental class of methodologies for time-dependent problems in fluid mechanics. Despite all favourable properties, these methods are not suitable for the computation of periodic flows. Adopting an existing space-time method does not ensure strict periodicity and demands excessive computational resources to establish a near-periodic signal via lengthy transient computations. In this paper, we circumvent these issues and propose to exploit the periodic nature by adopting a periodic space-time setting of arbitrary continuity via isogeometric analysis. Furthermore, we combine the usage of the RBVMS methodology and weak boundary conditions which provides a powerful space-time computational methodology. We show conservation properties of the proposed method and present a conservative traction evaluation. We use our computational setup for the simulation of incompressible flow past a prescribed periodically moving hydrofoil.

The paper is organized as follows. We present the time-periodic continuous space-time methodology in \cref{sec:domain}. In particular, we discuss the extraction of the mesh-constraint boundary velocity. Next, in \cref{sec:consprop} we discuss the conservation properties and present a force extraction method. In \cref{sec:results} we provide results of numerical experiments considering the mesh-constraint boundary velocity, force extraction, and periodic flow. We close with concluding remarks in \cref{sec:conclusions}.

\section{Periodic space-time formulation of the incompressible flow equations}
\label{sec:domain}

\subsection{Governing equations}
\label{sec:governingequations}
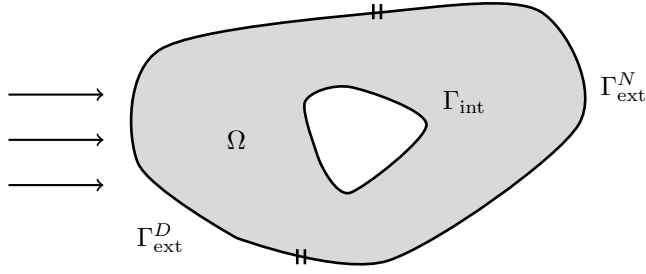
\begin{figure}[h!]
    \centering
\begin{tikzpicture}
\draw[black, line width=0.35mm, fill=gray!30] plot[smooth] coordinates {(3,1.5) (5,1.2) (7.5,3) (7,4.5) (5,4.5) (2,4) (1.7,2.5) (3,1.5)};
\draw[black, line width=0.35mm, fill=white] plot[smooth] coordinates {(4.1,2.5) (4.5,2.1)  (5.5,3.0)  (4.5,3.5)  (3.9,3.3) (4.1,2.5)};
\draw [black, line width=0.35mm] (3.8,1.35) -- (3.8,1.15);
\draw [black, line width=0.35mm] (3.9,1.35) -- (3.9,1.15);
\draw [black, line width=0.35mm] (4.8,4.4) -- (4.8,4.6);
\draw [black, line width=0.35mm] (4.9,4.4) -- (4.9,4.6);
\draw [->, thick](0, 2.2) -- (1.25, 2.2);
\draw [->, thick](0, 2.8) -- (1.25, 2.8);
\draw [->, thick](0, 3.4) -- (1.25, 3.4);
\node at (6, 3.3) {$\nsb$};
\node at (8.1, 3.5) {$\nb$};
\node at (2, 1.5) {$\db$};
\node at (3, 2.8) {$\Omega$};
\end{tikzpicture}
\caption{Sketch of the spatial domain with its boundaries, with inflow on the left.}
\label{tikz:problem}
\end{figure}

Consider a time-dependent spatial domain $\Omega = \Omega(t) \subset \mathbb{R}^d$ with boundary $\Gamma = \nsb \cup \Gamma_{\textrm{ext}}$ composed of a time-dependent interior $\nsb=\nsb(t)$ and exterior part $\Gamma_{\textrm{ext}}$. 
The outward unit normal to the boundary $\Gamma$ is defined as $\mathbf{n}$. Let us now consider a velocity field $\mathbf{u}$ and introduce the normal velocity $u_n = \mathbf{u}\cdot \mathbf{n}$ with positive and negative parts $u_n^{\pm}=\frac{1}{2}(u_n \pm |u_n|)$. 
We partition the exterior boundary into an inflow and outflow part according to the definitions:
\begin{subequations}
\begin{align}
    \Gamma_{\rm ext}^D :=&~ \left\{ \bx \in \Gamma|\un(\bx) < 0 \right\},\\
    \Gamma_{\rm ext}^N :=&~ \left\{ \bx \in \Gamma|\un(\bx) \geq 0 \right\}.
\end{align}
\end{subequations}
The domain is depicted in \cref{tikz:problem}.

We now consider the problem that reads in strong form:
\begin{subequations}\label{eq: strong form problem}
\begin{alignat}{2}
	\partial_t \bu + \bu \cdot \nabla  \bu + \nabla \p - \nabla \cdot \left(2 \nu \nabla^s \bu\right) 
	&= 
	\bfo 
	&& \quad \textrm{in} \quad \Omega 
	, \label{eq: strong form problem: lin mom} \\ 
	\nabla \cdot \bu 
	&= 
	0
	&& \quad \textrm{in} \quad \Omega 
	, \label{eq: strong form problem: cont}\\
	\bu 
	&= 
	\bg_{\mathrm{int}}
	&&\quad \textrm{in} \quad \nsb
	\label{eq: strong form problem: int}
	, \\
	\bu 
	&= 
	\bg_{\mathrm{ext}}
	&& \quad \textrm{in} \quad \db
	\label{eq: strong form problem: inflow}
	, \\
	-p \bn +\nu \nabla \bu \cdot \bn + \un^- \bu
	&= \bO
	&& \quad \textrm{in} \quad \nb
	\label{eq: strong form problem: outflow}
	, \\
	\bu(\cdot,0) &= \bu_0
    &&\quad \textrm{in} \quad \Omega.
    \label{eq: strong form problem: ic u}%
\end{alignat}
\label{eq:navstoinit}%
\end{subequations}
Here the unknown fields are the velocity $\mathbf{u}=\mathbf{u}(\mathbf{x},t) $ and the pressure $p=p(\mathbf{x},t) $ with spatial coordinate $\mathbf{x}$ and the time coordinate $t \in \mathcal{I}=(0,T)$ with final time $T>0$. We employ the standard notation for the gradient ($\nabla$), the symmetric gradient ($\nabla^s$) and the divergence ($\nabla\cdot$). Furthermore, $\nu$ denotes the (constant) kinematic viscosity, $\mathbf{f}=\mathbf{f}(t)$ is a (time-dependent) external force, and  $\mathbf{g}_{\textrm{int}}=\mathbf{g}_{\textrm{int}}(t)$ and $\mathbf{g}_{\textrm{ext}}$ are prescribed (time-dependent) velocities on the interior boundary and inflow partition of the exterior boundary, respectively. 
We split the prescribed no-slip velocity into a normal ($\bg_{n}$) and tangential component  ($\bg_{t}$):
\begin{subequations}
\begin{align}
    \bg_\textrm{int} &= \bg_n + \bg_{t},\\
    \bg_{n}  &=  (\bg_\textrm{int} \cdot\bn) \bn, \\
    \bg_{t} \cdot \bn &= 0.
\end{align}
\end{subequations}
Denoting the normal velocity of the domain boundary $\nsb$ by $v_n=\bg_n\cdot \mathbf{n}$, the normal component $\bg_{n}$ is prescribed by the relation $\bg_{n} = v_n \mathbf{n}$.

The equations \eqref{eq: strong form problem} describe the incompressible Navier-Stokes equations, with the balance of linear momentum and the continuity equation in \eqref{eq: strong form problem: lin mom} and \eqref{eq: strong form problem: cont}, the Dirichlet boundary conditions on the interior and the inflow boundary in \eqref{eq: strong form problem: int} and \eqref{eq: strong form problem: inflow}, the outflow boundary condition in \eqref{eq: strong form problem: outflow} and the initial condition in \eqref{eq: strong form problem: ic u}.

\subsection{Space-time formulation}\label{sec:st_form}
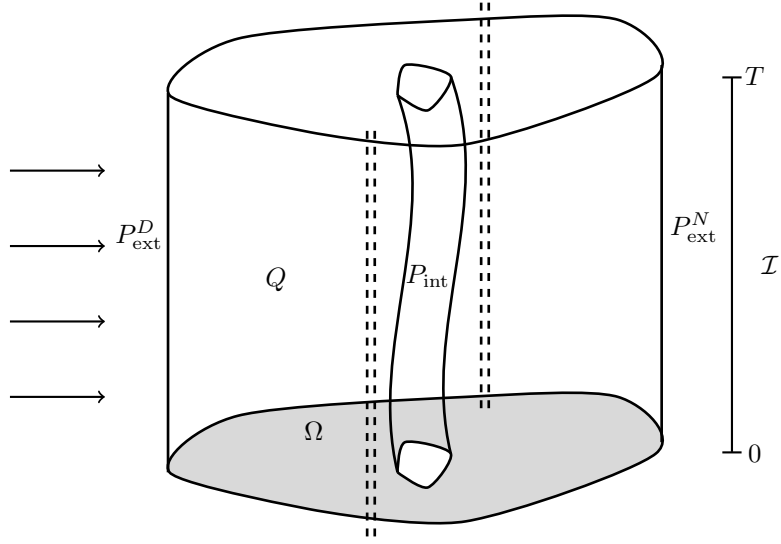
\begin{figure}[!ht]
    \centering
\begin{tikzpicture}
% Bot
\draw[black, line width=0.35mm, fill=gray!30] plot[smooth] coordinates {(3,0.75) (5,0.6) (7.5,1.5) (7,2.25) (5,2.25) (2,2) (1.1,1.25) (3,0.75)};
\draw[black, line width=0.35mm, fill=white] plot[smooth] coordinates {(4.1,1.25) (4.5,1.05) (4.8,1.5)  (4.2,1.65) (4.1,1.25)};
%Connection
\draw [black, line width=0.35mm] (1.08,1.25) -- (1.08,6.28);
\draw [black, line width=0.35mm](4.1,1.25) .. controls (3.7, 2.92) and (4.7, 4.6) .. (4.1,6.25);
\draw [black, line width=0.35mm](4.8,1.5) .. controls (4.4, 3.17) and (5.4, 4.83) .. (4.8,6.5);
\draw [black, line width=0.35mm] (7.575,1.65) -- (7.575,6.65);

%Top
\draw[black, line width=0.35mm] plot[smooth] coordinates {(3,5.75) (5,5.6) (7.5,6.5) (7,7.25) (5,7.25) (2,7) (1.1,6.25) (3,5.75)};
\draw[black, line width=0.35mm, fill=white] plot[smooth] coordinates {(4.1,6.25) (4.5,6.05) (4.8,6.5)  (4.2,6.65) (4.1,6.25)};
% Boundary separators
\draw [black, dashed, line width=0.35mm] (3.7,0.4) -- (3.7,5.85);
\draw [black, dashed, line width=0.35mm] (3.8,0.4) -- (3.8,5.85);
\draw [black, dashed, line width=0.35mm] (5.2,2.1) -- (5.2,7.5);
\draw [black, dashed, line width=0.35mm] (5.3,2.1) -- (5.3,7.5);
% Inflow arrows
\draw [->, thick](-1, 2.25) -- (0.25, 2.25);
\draw [->, thick](-1, 3.25) -- (0.25, 3.25);
\draw [->, thick](-1, 4.25) -- (0.25, 4.25);
\draw [->, thick](-1, 5.25) -- (0.25, 5.25);
% Time arrow
\draw [|-|, thick](8.5, 1.5) -- (8.5, 6.5);
% Labels
\node at (3, 1.8) {$\Omega$};
\node at (2.5, 3.8) {$Q$};
\node at (0.7, 4.4) {$\dbs$};
\node at (8.0, 4.5) {$\nbs$};
\node at (4.50, 3.85) {$\nsbs$};
\node at (9.0, 4.0) {$\mathcal{I}$};
\node at (8.8, 1.5) {$0$};
\node at (8.8, 6.5) {$T$};
\end{tikzpicture}
\caption{Sketch of the space-time domain $Q$ with its boundaries $P$, with inflow on the left, as an extrusion of the spatial domain $\Omega$ in gray.}
\label{tikz:stdom}
\end{figure}

We introduce the (continuous) space-time domain $Q = \Omega \times \mathcal{I}$ as an extrusion of the spatial domain $\Omega=\Omega(t)$. The boundary of $Q$ consists of an interior part $P_{\textrm int} =P_{\textrm int}(t) = \Gamma_{\rm int}(t) \times \mathcal{I}$, and an exterior part made up of an inflow $P_{\textrm ext}^D = \Gamma_{\rm ext}^D \times \mathcal{I}$ and an outflow $P_{\rm ext}^N = \Gamma_{\rm ext}^N \times \mathcal{I}$ contribution. We visualize the setup in \cref{tikz:stdom}.

We introduce the space-time coordinate $\bhatx= [\bx^T~\vst  t]^T = [x_1~ ...~ x_{d} ~\vst x_{d+1}]$ and the extended velocity vector $\bhatu = [\bu^T ~\vst]^T$, where $\vst$ is a velocity relating the time and space dimensions. For simplicity, $\vst$ can be chosen as $1$. 

In this work we focus on periodic flows and as such, we consider a periodically changing domain $\Omega$ with period $\mathcal{T}$:
\begin{equation}
    \Omega|_{t} = \Omega|_{ t + \mathcal{T}}.
\end{equation}
Additionally, we require the prescribed 
external force $\bfo$ and boundary velocities to be periodic:
\begin{subequations}
\begin{align}
\bfo (\bx, t) &= \bfo (\bx, t + \mathcal{T}),\\
\bg (\bx, t) &= \bg (\bx, t + \mathcal{T}).
\end{align}
\end{subequations}
The initial condition in \eqref{eq: strong form problem: ic u} is represented in the space-time setting by the time-periodic condition:
\begin{equation}
    \bu(\cdot,0) = \bu(\cdot,\mathcal{T}) \quad \textrm{in} \quad \Omega .
\end{equation}
We take the final time as $T = \mathcal{T}$ to cover one period.

Using these definitions, problem \eqref{eq: strong form problem} transforms in the space-time context into the steady state problem:
\begin{subequations}
\begin{alignat}{2}
    \bhatu \cdot \nabla_{\bhatx}  \bu
    + \nabla \p - \nu \nabla^2 \bu 
	&= 
	\bfo 
	&& \quad \textrm{in} \quad Q
	,  \label{eq:navstofinalmom}\\ 
	\nabla \cdot \bu 
	&= 
	0
	&& \quad \textrm{in} \quad Q
	, \\
	\bu 
	&= 
	\bg_{\mathrm{int}}
	&&\quad \textrm{in} \quad \nsbs
	\label{eq:nsb}
	, \\
	\bu 
	&= 
	\bg_{\mathrm{ext}}
	&& \quad \textrm{in} \quad \dbs
	\label{eq:ifb}
	, \\
	-p \bn +\nu \nabla \bu \cdot \bn + \un^- \bu
	&= \bO
	&& \quad \textrm{in} \quad \nbs
	,  \label{eq:navstofinalout}\\ 
	\bu(\cdot,0) &= \bu(\cdot,T) 
    &&\quad \textrm{in} \quad \Omega 
    .
\end{alignat}
\label{eq:navstostfinal}%
\end{subequations}
The first term in \cref{eq:navstofinalmom} represents the material derivative via
\begin{equation}
    \partial_t \bu + \bu \cdot \nabla  \bu= \bhatu \cdot \nabla_{\bhatx}  \bu,
    \label{eq:time-to-conv}
\end{equation}
where  $\nabla_{\bhatx}$ is the space-time gradient. The normal $\bn$ in \cref{eq:navstofinalout} is the classical spatial normal and can be extracted from the space-time normal $\bhatn = [n_1 ~... ~n_{d}~ n_{d+1}]^T$ via,
\begin{equation}
    \bn = \frac{1}{\sqrt{n_1^2+...+n_{d}^2}}\bmqty{n_1 \\ \vdots \\ n_d }. 
\end{equation}
The space-time outward normal $\bhatn$ has unit length in the norm $\|\cdot\|_{G_\vst}$ defined by
\begin{equation}
    \|\bhatn\|_{G_\vst}^2 = \bhatn \cdot \bG_\vst\bhatn,
\end{equation}
where $\bG_\vst$ is the space-time metric
\begin{equation}
    \bG_\vst = 
    \left (\begin{array}{cc}
{\bf{I}}_{d \times d}& 0_{1\times d} \\
0_{d \times 1} & \vst^2 
\end{array}\right ).
\end{equation}
Furthermore, the normal velocity $v_n$ is related to the space-time velocity $s$ and the space-time normal $\hat{\mathbf{n}}$ via:
\begin{equation}\label{eq:vn}
v_{n} = -\vst\frac{n_{d+1}}{\sqrt{n^2_{1} + .. + n^2_{d}}}.
\end{equation}

\subsection{Weak formulation of the continuous space-time problem}
The weak formulation of the continuous space-time problem is stated using the trial and test function spaces $\WW_g$ and $\WW_0$ respectively. 
Members of the trial function space $\WW_g$ satisfy the non-homogeneous Dirichlet boundary conditions for the velocity on $\dbs$ whereas
elements in the test function space $\WW_0$ satisfy the homogeneous Dirichlet boundary conditions on $\dbs$. Additionally, members of both spaces satisfy the periodic boundary condition $\bu |_{\Omega_0} = \bu |_{\Omega_T}$ where,
$\Omega_0 = Q|_{t=0}$ and $\Omega_T = Q |_{t=T}$. To enforce the Dirichlet boundary conditions on $\nsbs$, we introduce the subspaces $\mathcal{V}_g \subset \WW_g$ and $\mathcal{V}_0 \subset \WW_0$ that additionally satisfy non-homogeneous and homogeneous boundary conditions on $\nsbs$, respectively.

The variational formulation of \cref{eq:navstostfinal} now reads as:\\ 

\textit{find  $ \bU = \left\{\bu, p\right\} \in \mathcal{V}_g $ such that for all  $ \bW = \left\{\bw, q \right\} \in   \mathcal{V}_0:$}

\begin{subequations}
\begin{equation}
	B_{\mathrm{GAL}}\left(\bU, \bW \right)
	=
	L \left(\bW\right),
\end{equation}
\\	
\indent \textit{where}
\\
\begin{align}
\begin{split}
B_{\mathrm{GAL}}\left(\bU, \bW \right),
=&  
	\left(\bw, \bhatu \cdot \nabla_{\bhatx}  \bu\right)_{Q} % time = convection notation
	+ 
	\left(\divbw, p\right)_{Q}
	\\&+
	\left(\gbw, \nu \nabla \bu\right)_{Q}
	+
	\left(q, \nabla \cdot \bu\right)_Q
	- 
    \left(\bw, \un^{-} \bu \right)_{\nbs}
\end{split}
\\[10pt]
L \left(\left\{\bw, q\right\}\right)
=& 
\left(\bw, \bfo\right)_{Q}.
\end{align}
\end{subequations}
The $L^2$ inner product over $D$ is defined as $(\cdot , \cdot)_D $.

\subsection{Weak formulation of the discrete problem}
\label{subsec: weak discrete}
To introduce the numerical discretization, we first subdivide our physical domain $Q$ into elements $Q_K$. The domain of element interiors denotes:
\begin{equation}
    \tilde{Q} = \bigcup_K Q_K. 
\end{equation}
We apply residual-based variational multiscale turbulence modeling \cite{Bazilevs2007,Bazilevs2010a} in which the weighting function space and trial solution space are decomposed into subspaces that contain the coarse and fine scales:
\begin{subequations}\label{eq: MS split}
\begin{alignat}{2}
    \WW_g =&~ \WW_g^h \oplus \WW',
    \\
    \WW_0 =&~ \WW_0^h \oplus \WW',
    \label{eq:smallcoarse}
\end{alignat}
\end{subequations}
where $\WW_g^h$ and $\WW^h_0$ are coarse-scale spaces, and $\WW' \subset \WW_g \cup \WW_0$ are the fine scales. The coarse-scale space is spanned by
the finite-dimensional numerical discretization whereas the fine-scales are their infinite-dimensional complement. Uniqueness of the multi-scale split \eqref{eq: MS split} is ensured when the split is established via a projection operator. \eqref{eq:smallcoarse} implies that the members of $\WW_g$ and $\WW_0$ split as:
\begin{subequations}
\begin{align}
    \left\{\bu, p\right\} &= \left\{\bu^h, p^h\right\} + \left\{\bu', p'\right\},
    \\
    \left\{\bw, q\right\} &= \left\{\bw^h, q^h\right\} + \left\{\bw', q'\right\}
    ,
\end{align}
\end{subequations}
where the components of the coarse-scale subspaces are denoted as $\bU^h = \left\{\bu^h, p^h\right\} \in \WW_g^h$ and  $\bW^h = \left\{\bw^h, q^h\right\} \in \WW^h_0$, and the components of the small-scale subspace are denoted as $\bU' = \left\{\bu', p'\right\} \in \WW'$ and  $\bW' = \left\{\bw', q'\right\} \in \WW'$.

To arrive at the fully-discrete formulation we make the following modeling choices. First, we apply a pseudo-transient continuation to march in pseudo-time to the space-time steady-state solution. Next, we select a standard $H_0^1$-multiscale projector that eliminates the fine-scale viscosity contribution. Next, we replace the small-scale space $\mathcal{W}'$ by the velocity-pressure product $\mathcal{V}'\times \mathcal{P}'$. The fine-scales are modeled as:
\begin{subequations}
\begin{alignat}{2}
	\bu' &= -\tau_M \br_M,\\
	p' &= -\tau_C r_C,
\end{alignat}
\end{subequations}
with the strong residuals
\begin{subequations}
\begin{alignat}{2}
	\br_M &= 	
	 \left(\bhatu^h \cdot \nabla_{\bhatx} \right) \bu^h 
	- \nabla \p - \nu \nabla^2 \bu^h - \bfo,  \\ 
	r_C &= \nabla \cdot \bu^h,
\end{alignat}
\end{subequations}
and stability parameters
\begin{subequations}
\begin{alignat}{2}
    \tau_M &= 
    \left(\bhatu^h \cdot \bhatG  \bhatu^h + C^I \nu^2 \bG:\bG\right)^{1/2}
    , \\
    \tau_C &= \tau_M^{-1} \mathrm{Tr}(\bG)^{-1}.
\end{alignat}
\end{subequations}
In both the momentum residual and its corresponding stability parameter the time derivative is incorporated in the convection term, analogous to \eqref{eq:time-to-conv}.
As a consequence, the convective and diffusive contributions depend on two different metric tensors, the space-time metric tensor $\bhatG$ and spatial metric tensor 
$\bG$, respectively. These metric tensors are given by
\begin{align}
    \bhatG = \left(\pdv{\bxi}{\bhatx}\right)^T \bG_\vst \pdv{\bxi}{\bhatx},
        &&
    \bG = \left(\pdv{\bxi}{\bx}\right)^T \pdv{\bxi}{\bx}.
\end{align}
Lastly, we enforce the Dirichlet boundary conditions weakly \cite{Bazilevs2007a}. To this purpose we introduce the penalty parameter $\tau_b$:
\begin{equation}
  \tau_b =  \frac{1}{2} C^{I}_b \nu \left(\bn \cdot \bG  \bn\right)^{\frac{1}{2}}. 
\end{equation}

We now define the following fully-discrete formulation:\\

\textit{find  $ \bU^h = \left\{\bu^h, p^h\right\} \in \WW_g^h $ such that for all  $ \bW = \left\{\bw^h, q^h \right\} \in   \WW^{h}_0$:}

\begin{subequations}\label{eq: discrete weak formulation}
\begin{equation}
	B\left(\bU^h, \bW^h \right)
	=
	L \left(\bW^h\right),
\end{equation}
\\
\indent \textit{where}
\\
\begin{align}
\begin{split}
B\left(\bU^h, \bW^h\right)
=& 
B_{\mathrm{GAL}}\left(\bU^h, \bW^h \right)
+
B_{\mathrm{PT}}\left(\bU^h, \bW^h \right) 
\\ &
+
B_{\mathrm{STAB}}\left(\bU^h, \bW^h \right)
+ 
B_{\mathrm{WBC}}\left(\bU^h, \bW^h \right),
\end{split}
\\[15pt]
\begin{split}
B_{\mathrm{PT}}\left(\bU^h, \bW^h \right) 
=&
\left(\bw^h, \partial_{\theta} \bu^h \right)_{Q}
+\frac{1}{a^2}\left(q^h, \partial_{\theta} \p^h \right)_{Q},
\label{eq:discweak:pt}
\end{split}
\\[10pt]
\begin{split}
B_{\mathrm{STAB}}\left(\bU^h, \bW^h \right)
=&
    -
    \left( \gbwhat^h, \sbu \otimes \bhatu^h \right)_{\tilde{Q}} % SUPG
    - 
    \left(\nabla \bw^h,  \bu^h \otimes \sbu \right)_{\tilde{Q}}  % Cross term
    \\&- 
    \left(\gbw^h,  \sbu \otimes \sbu\right)_{\tilde{Q}}
    -
    \left(\gbq^h, \sbu \right)_{\tilde{Q}} % PSPG
    - 
    \left(\divbw^h,  p'\right)_{\tilde{Q}} % LSIC
    ,
\label{eq:discweak:stab}
\end{split}
\\[10pt]
\begin{split}
B_{\mathrm{WBC}}\left(\bU^h, \bW^h \right) 
=&
\left(\bw^h, p^h \bn - 
\nu \nabla_{\bx} \bu^h \cdot \bn
\right)_{\nsbs} % Consistency
+
\left(
\nu \nabla \bw^h \cdot \bn
- q^h\bn, \bu^h -\bg \right)_{\nsbs} % Dual Consistency
\\&+
\left(\bw^h \tau_b,  \bu^h -\bg\right)_{\nsbs}. % Penalty
\label{eq:discweak:wbc}
\end{split}
\end{align}
\label{eq:discweak}
\end{subequations}

Equation \eqref{eq:discweak:pt} represents the pseudo-transient continuation as a globalization technique \cite{Kelley1998,Coffey2003}. Pseudo-transient continuation technique is a widely applied methodology that obtains the steady-state solution by adding a derivative to pseudo-time $\theta$. The first term is classical, whereas the utilization of the second term is non-standard. This term introduces artificial compressibility \cite{Chorin1997,Chorin1968,Temam1969}, 
where $a$ is an artificial speed of sound. This term overcomes some of the difficulties due to the saddle-point nature of the underlying problem (i.e. the absence of a pressure term in the continuity equation). Moreover, we note the introduction of this term permits more powerful preconditioning options such as algebraic multigrid (AMG). We remark that the numerical solution of the problem is fully incompressible and thus does not depend on the artificial speed of sound $a$. 

Equation \cref{eq:discweak:stab} describes terms associated with variational multiscale stabilisation \cite{Bazilevs2007}. 
In LES terminology the first two terms represent the cross-stress, while the third term represents the Reynolds stress. In the context of stabilized methods, the first term is the Streamline-upwind Petrov-Galerkin (SUPG) term, and the fourth and last terms are the PSPG and LSIC terms respectively. Note that the first and the second terms are not each other transposes. Namely, we incorporate the temporal derivative of the fine-scales in the SUPG term:
\begin{equation}
    (\bw, \partial_t \sbu)_{\tilde{Q}} 
    + 
    \left( \gbw^h, \sbu \otimes \bu^h \right)_{\tilde{Q}}
    =
    \left( \gbwhat^h, \sbu \otimes \bhatu^h \right)_{\tilde{Q}}.
\end{equation}
This relation is a direct consequence of the apply partial integration (in the temporal direction) of the fine-scale time-derivative term:
\begin{equation}
    (\bw, \partial_t \sbu)_{\tilde{Q}} 
    = 
    -(\partial_t \bw, \sbu)_{\tilde{Q}},
    \label{eq:pismallscales}
\end{equation}
where we note the absence of boundary contributions due to the periodic boundary conditions.

Lastly, equation \cref{eq:discweak:wbc} enforces the weak boundary conditions on the interior boundary \eqref{eq:nsb}. The first term is the consistency term. This term originates from integration by parts and as such guarantees variational consistency. The second term is the so-called the dual consistency term, and the last term is the penalty term that ensures the stability of the formulation. We recall that the Dirichlet boundary conditions in \eqref{eq:ifb} on $\dbs$  are enforced strongly.

\section{Conservation properties}
\label{sec:consprop}
In this section we establish the conservation properties of the discrete method. We show conservation of mass, conservation of linear momentum and provide an approach to conservatively evaluate the traction. We consider a converged solution where $\partial_\theta \bu^h = \partial_\theta p^h =0$.

\subsection{Conservation of mass}
The global conservation of mass directly follows by selecting the weighting function $\bW^h = \left\{\mathbf{0},1\right\}$ in the discrete weak formulation \eqref{eq: discrete weak formulation}:
\begin{align}
  \displaystyle\int_Q \nabla \cdot \bu^h ~{\rm d}x = 0.
\end{align}
We do not attain conservation of mass per time-slab since the weighting function with pressure component that equals $1$ on a single time-slab and $0$ on the others is not a member of $\WW_0^h$. Remark that it is possible to work with a particular selection of isogeometric velocity-pressure spaces that establishes pointwise satisfaction of the incompressibility constraint \cite{evans2013isogeometric,EiAk17ii}.

\subsection{Conservation of linear momentum}

In order to study the conservation of linear momentum one might wish to substitute the weighting function $\bW^h = \left\{\bw^h,q^h\right\} = \left\{ \mathbf{e}_i, 0\right\}$ with $\mathbf{e}_i \in \mathbb{R}^{d}$ the $i$-th Cartesian unit vector into the discrete weak formulation \eqref{eq: discrete weak formulation}. This choice is not permitted: $\left\{ \mathbf{e}_i, 0\right\} \notin \WW^h_0$. One possible remedy is to work with unconstrained function spaces and weakly enforce the non-homogeneous boundary condition via a Lagrange multiplier construct \cite{EiAk17i,EiAk17ii}. The Lagrange multiplier is also called \textit{auxiliary flux} \cite{Hughes2000} and is used to show global and local conservation. The method yields conservative boundary fluxes which is a major advantage as compared to utilizing direct procedures that provide non-conservative boundary fluxes. 

We denote the vector-valued Lagrange multiplier/auxiliary flux as $\boldsymbol{\lambda}$. Recall that the discrete weak formulation \eqref{eq: discrete weak formulation} is defined for the test function space $\WW^{h}_0$ in which the velocity test functions vanish on $\dbs$. In order to present the augmented formulation, we require the introduction of other test function spaces. Denote the set of all velocity basis functions $\eta$ and furthermore denote with $\eta_g$ the set of velocity basis functions that do not vanish on $\dbs$. With the notation 
$\WW^{h}_0=\UU^{h}_0\times \PP^h$ of the velocity and pressure components of the test function space, we have $\UU^h_0 = {\rm span}\left\{\boldsymbol{N}_A\right\}_{A \in \eta - \eta_g}$, where $\boldsymbol{N}_A=\boldsymbol{N}_A(\boldsymbol{x})$ are the velocity basis functions. Furthermore, we introduce the unrestricted velocity space $\UU^h = {\rm span}\left\{\boldsymbol{N}_A\right\}_{A \in \eta }$ and unrestricted velocity-pressure space $\WW^h = \UU^{h}\times \PP^h$. The augmented problem now reads:\\

\textit{find  $ \bU^h \in \WW_g^h $ such that for all  $\hat{\bW}^h = \left\{\hat{\bw}^h,q^h\right\} \in \WW^{h}$:}
\begin{align}\label{eq: augmented problem}
  (\boldsymbol{\lambda}^h,\hat{\bw}^h)_{\dbs} =	B\left(\bU^h, \hat{\bW}^h \right) - L \left(\hat{\bW}^h\right).
\end{align}
This problem splits as:\\

\textit{find  $ \bU^h \in \WW_g^h $ and $\boldsymbol{\lambda}^h \in \WW^{h}-\WW^h_0$ such that}
\begin{subequations}
\begin{align}
  0 =&~	B\left(\bU^h, \bW^h \right) - L \left(\bW^h\right) \quad  \text{for all} \quad \bW^h \in \WW^{h}_0\\
  (\boldsymbol{\lambda}^h,\hat{\bw}^h)_{\dbs} =&~	B\left(\bU^h, \hat{\bW}^h \right) - L \left(\hat{\bW}^h\right) \quad  \text{for all} \quad \hat{\bW}^h \in \WW^{h}-\WW^h_0.
\end{align}
\end{subequations}
The first subproblem coincides with our original weak formulation and thus completely determines the numerical solution $\bU^h \in \WW_g^h$. This solution may be directly substituted into the second subproblem to evaluate the discrete auxiliary flux $\boldsymbol{\lambda}^h \in \WW^{h}-\WW^h_0$.

We are now in the position to evaluate the linear momentum conservation and select $\hat{\bW}^h = \left\{ \mathbf{e}_i, 0\right\}$ in \eqref{eq: augmented problem}:
\begin{align}\label{eq: lin mom}
  \displaystyle\int_{\dbs} \lambda^h_i ~{\rm d}s =&~	\displaystyle\int_{\nsbs} p^h n_i - 
\nu (u_{i,j}+u_{j,i})n_j
~{\rm d}s - \displaystyle\int_Q f_i ~{\rm d}x ~{\rm d}s\nonumber\\
&~-
\displaystyle\int_{\nbs}\un^- u_i^h~{\rm d}s+ \displaystyle\int_{\nsbs} \tau_b (u_i^h -g_i)~{\rm d}s.
\end{align}
This shows that $\lambda^h_i$ represents the total conserved boundary flux on $\dbs$. Remark that the latter two members on the right-hand side result from the usage of weak boundary conditions on $\nsbs$ and are thus absent when instead imposing these conditions strongly. 

\subsection{Conservative traction evaluation}
\label{subsec:Conservative traction evaluation}
With the aim of evaluating the time-dependent traction on the interior boundary $\nsb$ we select $\hat{\bW}^h = \left\{ \mathbf{e}_i N_a, 0\right\}$ in \eqref{eq: augmented problem} with $N_a=N_a(z)$ an arbitrary basis function in the temporal direction. Note that this choice is permitted due to the tensor structure of the NURBS computational mesh. Substitution provides:
\begin{align}\label{eq: traction1}
  \displaystyle\int_{\dbs} \lambda^h_i N_a ~{\rm d}s + \displaystyle\int_Q f_iN_a ~{\rm d}x+\displaystyle\int_{\nbs}\un^- u_i^h~{\rm d}s=&~
  \displaystyle\int_{\nsbs} p^h n_i N_a - 
\nu (u_{i,j}+u_{j,i})n_j N_a
~{\rm d}s \nonumber\\
&~+ \displaystyle\int_{\nsbs} \tau_b (u_i^h -g_i)N_a~{\rm d}s.
\end{align}

The right-hand side of \eqref{eq: traction1} contains all the integrals on the interior boundary $\nsbs$. In order to evaluate the (vector-valued) traction force $\boldsymbol{\psi}$ we introduce the discrete problems for $i=1,~ \hdots,~ d$:

\textit{find  $\psi^h_i \in {\rm span}\left\{N_b\right\}_{b \in \xi}$ such that}
\begin{align}\label{eq: traction2}
  \displaystyle\int_{\nsbs} \psi^h_i N_a ~{\rm d}s =&~
  \displaystyle\int_{\nsbs} p^h n_i N_a - 
\nu (u_{i,j}+u_{j,i})n_j N_a
~{\rm d}s \nonumber\\
&~+ \displaystyle\int_{\nsbs} \tau_b (u_i^h -g_i)N_a~{\rm d}s
\end{align}

where $\xi$ is the set of basis function numbers in the time direction. The traction forces $\psi^h_i$ thus result from inverting a mass matrix (per direction).

\section{Numerical experiments}
\label{sec:results}
In this section, we discuss the computational setup and subsequently provide results of four numerical experiments using the formulation in \cref{subsec: weak discrete}.  We evaluate the forces in the space-time domain using the conservative traction evaluation of \cref{subsec:Conservative traction evaluation}. First, we compare the results of the mesh-constraint boundary velocity of a sinusoidal heaving hydrofoil with the analytical solution and study its dependency on the temporal discretization. Second, in order to examine the capability of the proposed methodology of predicting steady flow, we study the results of fluid flow past a stationary hydrofoil. We perform a grid convergence study and compare our results with the literature. Third, we focus on the hydrodynamics of a moving body, which is much more complex than the case of a steady body. We simulate the flow past a low-frequency heaving hydrofoil. Lastly, we investigate the predictive capability of the methodology on capturing history effects in the wake. We simulate the flow past a pitching hydrofoil at a moderate frequency. Experimental data considering (unsteady) forces on a hydrofoil in a low Reynolds-number flow is not available in the literature. We support our predictions with numerical results from the literature and two-dimensional simulations of the steady variant of the flow model \eqref{eq:navstoinit} using the same spatial discretization.

\subsection{Computational setup}
\begin{figure}[h]
\centering
\begin{tikzpicture}
% Define coordinates
\coordinate (A) at (-0.46,2);
\coordinate (B) at (3.0,2);
\coordinate (C) at (3.0, -2);
\coordinate (D) at (-0.46,-2);
\coordinate (le) at (-0.66,0);
\coordinate (te) at (0.66,0);
% Draw domain and boundaries
\filldraw[black, line width=0.35mm, fill=gray!30]
% Outer boundary
(A) -- (B) -- (C) -- (D) .. controls (-4,-2) and (-4,2) .. cycle
% Inner boundary
(te) .. controls (-0.33,  0.133) and  (-0.66,  0.133) .. (le) .. controls (-0.66, -0.133) and (-0.33, -0.133)  .. cycle ;
% Draw annotations
\node at (0, 2.3) {$\db$};
\node at (-2.8, -1.8) {$\db$};
\node at (-3.2, 2.2) {$\nsb$};
\draw [-, thin](-2.9, 1.9) -- (-0.7, 0.2);
\node at (3.5, 0.0) {$\nb$};
\node at (1.5, 0.5) {$\Omega$};
% Inflow arrows
\draw [->, thick](-5, 1.0) -- (-3.75, 1.0);
\draw [->, thick](-5, 0.0) -- (-3.75, 0.0);
\draw [->, thick](-5, -1.0) -- (-3.75, -1.0);
% Nurbs patches boundaries
\draw[black, line width=0.35mm, dotted] (te) -- (3.0, 0,0);
\draw[black, line width=0.35mm, dotted] (te) -- (0.66, -2,0);
\draw[black, line width=0.35mm, dotted] (te) -- (0.66, 2,0);
\draw[black, line width=0.35mm, dotted] (le) -- (-3.1, 0,0);
\draw[black, line width=0.35mm, dotted] (-0.46, -0.12) -- (-0.46, -2.0);
\draw[black, line width=0.35mm, dotted] (-0.46, 0.12) -- (-0.46, 2.0);
\end{tikzpicture}
\caption{Schematic representation of the domain $\Omega$, as a time slice of $P$, surrounding the hydrofoil with the no-slip boundary $\nsb$, the inflow boundary $\db$ and the outflow boundary $\nb$. The arrows indicate the direction of the flow. The six NURBS patches are indicated with a dotted line.}
\label{tikz:problemfoil}
\end{figure}
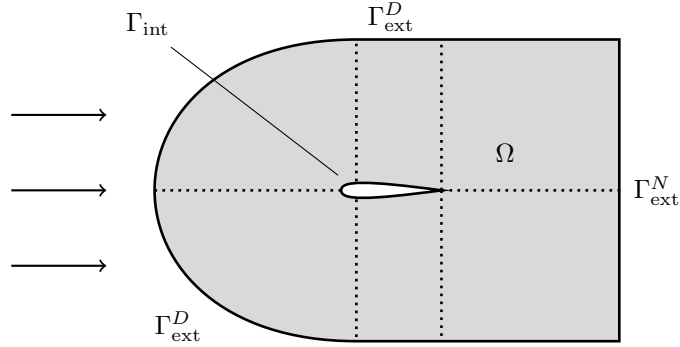

We introduce the space-time domain $Q$ as an extrusion of the spatial domain $\Omega$ enclosing a symmetric four-digit NACA foil section \cite{Abbott1959}. The spatial domain is discretized as a C-shaped mesh using six NURBS patches employing second-order NURBS. The spatial domain is illustrated in \cref{tikz:problemfoil}. The discretization is $\mathrm{C}^1$-continuity inside the patches and $\mathrm{C}^0$-continuity across patches. The hydrofoil and its motion are incorporated into the space-time mesh using curve interpolation.

\begin{figure}[h]
\centering
\begin{subfigure}{0.49\textwidth}
\centering
\includegraphics[width=0.83\textwidth, trim={5cm 7.0cm 2cm 6.0cm}, clip]{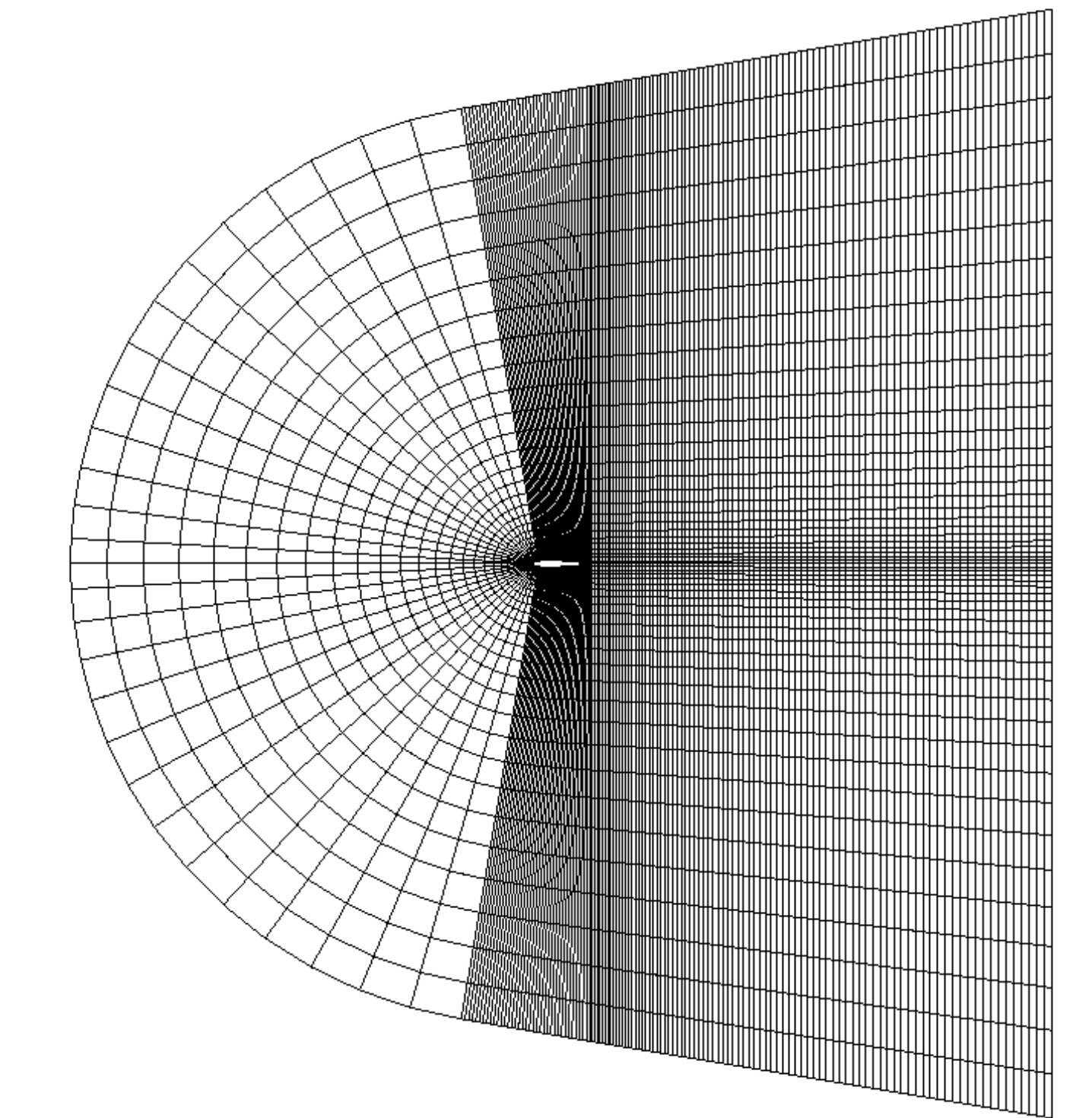}
\caption{}
\label{fig:mesh}
\end{subfigure}%
\hfill
\begin{subfigure}{0.49\textwidth}
\centering
\includegraphics[width=\textwidth, trim={5cm 7.5cm 0cm 6.5cm}, clip]{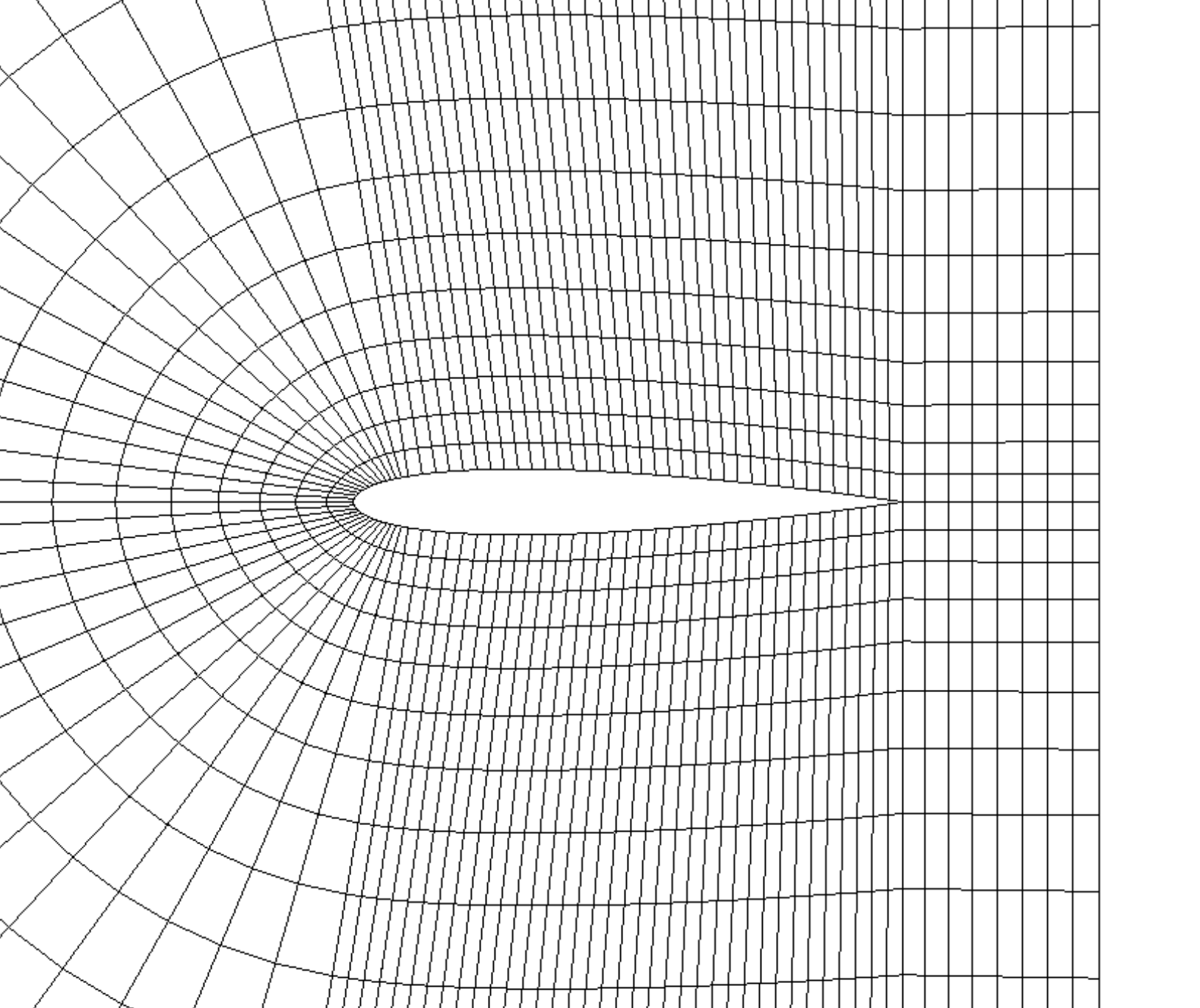}
\caption{}
\label{fig:meshzoom}
\end{subfigure}
\caption{The spatial mesh as a slice in time: a) The full C-shaped spatial mesh; b) A close-up of the spatial mesh near the interpolated hydrofoil.}
\label{fig:meshes}
\end{figure}

\cref{fig:meshes} provides an overview and a close-up of a temporal slice of the mesh. The mesh is constructed with the aim of achieving high quality near the hydrofoil. Based on simulations of the flow past a cylinder \cite{Behr1995}, we choose the distance between $\nsb$ and $\eb$ to be 8 chord lengths in order to preclude influence from the outflow boundary $\eb$. We have numerically verified that influence of $\eb$ is virtually absent. We select the chord $c$ and free stream velocity $U$ as $c = U = 1$. The numerical experiments are conducted in DelFI, which is based on the MFEM library \cite{Anderson2021}. 

The time-marching in pseudo-time $\theta$ towards a steady solution typically consists of $14$ pseudo-time steps of $5$ seconds using the backward Euler method as a pseudo-time marching scheme. As a stopping criterion, we terminate the computation when the $L^2$-norm of the residual of the momentum and mass equations is smaller than $10^{-6}$ at the start of the first Newton iteration. Per time step we use 5 Newton iterations.
We choose the artificial speed of sound $a$ as $4$, which exceeds the velocities encountered in the simulations. This provides a significant reduction in simulation time. Furthermore, we select the inverse estimate coefficients as $C^I = 36$ and $C^{I}_b = 8$. We note that the latter is only suitable for polynomial degrees up to $2$.

\pgfplotsset{every axis/.append style={
    cycle list/Dark2,
    cycle multiindex* list={
            mark list*\nextlist
            Dark2\nextlist},
    line width=1.0pt,
    mark size=2.5pt,
    grid = both,
    minor tick num = 1,
    major grid style = {lightgray},
    minor grid style = {lightgray!25},
    ylabel style={yshift=-0.20cm},
    yticklabel style={
    /pgf/number format/precision=4,
    /pgf/number format/fixed},
    xticklabel style={
    /pgf/number format/precision=3,
    /pgf/number format/fixed},
    tick label style={font=\tiny},
    label style={font=\footnotesize},
    legend style={font=\scriptsize},
    legend cell align={left},
    }}

Lastly, we discuss the computation of the boundary velocity on the interior boundary $\bg_\mathrm{int}$. We recall the split:
\begin{subequations}
\begin{align}
    \bg_\textrm{int} &= \bg_n + \bg_{t},\\
    \bg_{n}  &=  (\bg_\textrm{int} \cdot\bn) \bn = v_n \bn, \\
    \bg_{t} \cdot \bn &= 0,
\end{align}
\end{subequations}
where $v_n$ satisfies the relation \cref{eq:vn}. The domain motion fully prescribes $\bg_{n}$, while the tangential component $\bg_{t}$ is still undetermined. To numerically determine $\bg_\textrm{int}$ however, we use the motion encoded in the mesh and do not rely on the relation \cref{eq:vn}. The following procedure is permitted due to the extrusion structure of the space-time mesh. We have the following relations:
\begin{subequations}
\begin{align}
t= \ &t(\xi_{d+1}),\\
\bX=\  &\bX(\xi_1,...\xi_d),
\end{align}
\end{subequations}
where $\bX$ is a Lagrangian coordinate labeling a particle, and where $\boldsymbol{\xi}$ are the coordinates in the reference domain. 
We compute the boundary velocity by taking the derivative of the spatial coordinate $\bx$ to the time direction $t = x_{d+1}$ on a particle path:
\begin{equation}
   \bg_\textrm{int}  = \left . \pdv{\bx}{t}\right |_{\bX} \quad \textrm{in} \quad \nsbs .
   \label{eq:boundvel}
\end{equation}
Realizing the dependence $\bx=\bx(\xi_1,...\xi_{d+1})$, we can use the chain rule to conclude:
\begin{equation}
   \bg_\textrm{int}  = \vst \left . \pdv{\bx}{\st}\right |_{\bX} 
   =  \vst \sum_{i=1}^{d+1} \pdv{\bx}{\xi_i} \left .\pdv{\xi_i}{\st} \right |_{\bX} 
   = \vst \pdv{\bx}{\xi_{d+1}} \pdv{\xi_{d+1}}{\st}  \quad \textrm{in} \quad \nsbs.
   \label{eq:boundvel_simple}
\end{equation}
We note that the velocity $\bg_\textrm{int}$ computed via \cref{eq:boundvel_simple} satisfies $\bg_\textrm{int} \cdot \bn= v_n$, where $v_n$ is given by \cref{eq:vn}.

\subsection{Mesh-constraint boundary velocity}

We evaluate the mesh motion and the resulting mesh-constraint boundary $\bg_\mathrm{int}$ velocity. We apply a heave motion to the hydrofoil such that it only moves in the $x_2$-direction. The heave motion of the hydrofoil is sinusoidal with $h(t) = h_a \sin (2\pi t/T)$, where the amplitude $h_a$ = 0.5 $m$ and the period $T$ = 8 $s$. We use three different temporal resolutions consisting of $6$, $12$ and $24$ elements in the temporal direction $n_{\mathrm{el}, x_3}$.

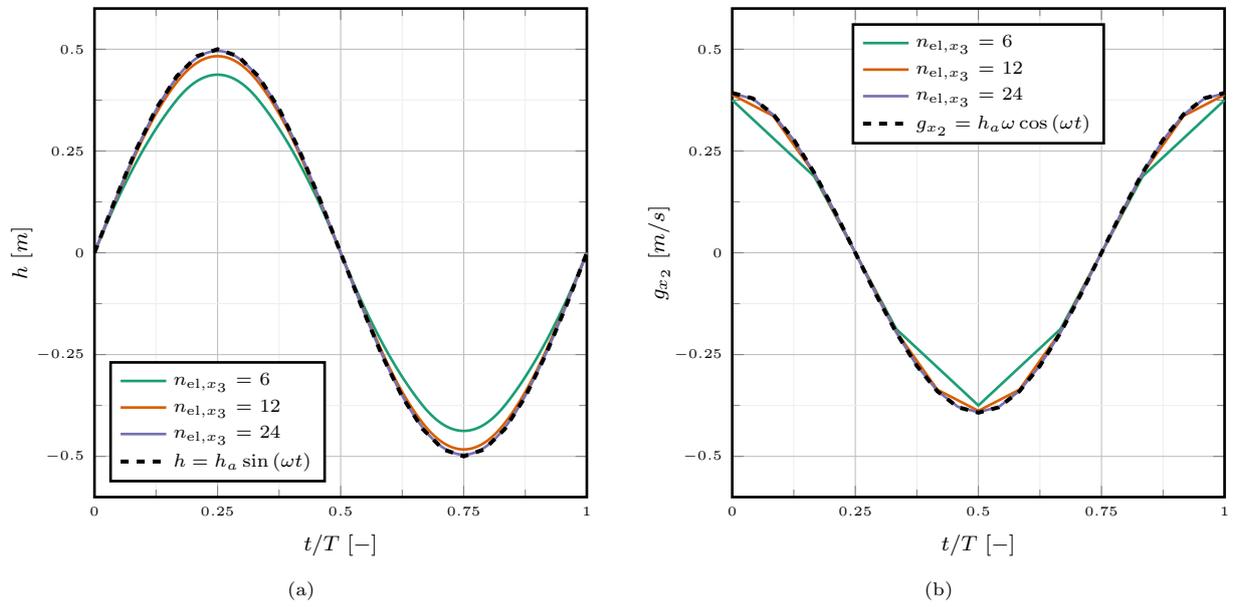
\begin{figure}[h]
    \centering
    \begin{subfigure}[b]{0.49\textwidth}
        \centering
        \begin{tikzpicture}
            \def\T{8.0}
            \def\ha{0.5}
            \def\omegasim{2*pi/\T}
            \begin{axis}[
            %only marks,
            xmin = 0, xmax = 1,
            ymin = -0.6, ymax = 0.6,
            xtick distance = 0.25,
            ytick distance = 0.25,
            width = \textwidth,
            height = \textwidth,
            xlabel = {$t/T$ $[-]$},
            ylabel = {$h$ $[m]$},
            domain=0:1,
            legend pos = south west,
            no markers,
            ]
            % x expr=\thisrow{x}/\T
            \addplot +[]table[x expr=\thisrow{x}/\T,y index=1] {data/bndfit/line-ncpz08.dat};
            \addplot +[]table[x expr=\thisrow{x}/\T,y index=1] {data/bndfit/line-ncpz14.dat};
            \addplot +[]table[x expr=\thisrow{x}/\T,y index=1] {data/bndfit/line-ncpz26.dat};
            \addplot [line width=1.5pt, mark=none, black, dashed,]{\ha*sin(deg(\omegasim*x*\T))};
            \legend{ $n_{\mathrm{el}, x_3}$ = 6 , $n_{\mathrm{el}, x_3}$ = 12, $n_{\mathrm{el}, x_3}$ = 24, $h = h_a   \sin\left(\omega t\right)$}
            \end{axis} 
        \end{tikzpicture}
    \caption{}
    \label{tikz:boundfit}
    \end{subfigure}
    \hfill
    \begin{subfigure}[b]{0.49\textwidth}
        \centering
        \begin{tikzpicture}
            \def\T{8}
            \def\ha{0.5}
            \def\omegasim{2*pi/\T}
            \begin{axis}[
            %only marks,
            xmin = 0, xmax = 1,
            ymin = -0.6, ymax = 0.6,
            xtick distance = 0.25,
            ytick distance = 0.25,
            width = \textwidth,
            height = \textwidth,
            xlabel = {$t/T$ $[-]$},
            ylabel = {$g_{x_2}$ $[m/s]$},
            domain=0:1,
            %legend pos = south west,
            legend style={at={(0.5,0.97)},anchor=north},
            no markers,
            ]
            \addplot +[]table[x expr=\thisrow{z}/\T, y index=1,y index=4] {data/bndvel/ncpz08.dat};
            \addplot +[]table[x expr=\thisrow{z}/\T,y index=4] {data/bndvel/ncpz14.dat};
            \addplot +[]table[x expr=\thisrow{z}/\T,y index=4] {data/bndvel/ncpz26.dat};
            \addplot [line width=1.5pt, mark=none, black, dashed,]{\ha*\omegasim*cos(deg(\omegasim*x*\T))};
            \legend{ $n_{\mathrm{el}, x_3}$ = 6 , $n_{\mathrm{el}, x_3}$ = 12, $n_{\mathrm{el}, x_3}$ = 24, $g_{x_2} = h_a \omega  \cos\left(\omega t\right)$}
            \end{axis}
        \end{tikzpicture}
    \caption{}
    \label{tikz:boundvelres}
    \end{subfigure}
    \caption{The motion $h$ in (a) and velocity $g_{x_2}$ in (b) of the hydrofoil in $x_2$-direction for 3 resolutions $n_{\mathrm{el}, x_3}$ in time direction and the analytical solution for a heave motion with $T$ = 8 $s$.}
    \label{fig:heavcurves}
\end{figure}

\cref{tikz:boundfit} presents the resulting mesh motion with the corresponding analytical solution. The second-order NURBS curves are reconstructed using the control points from the mesh. We observe that the finest mesh with $n_{\mathrm{el}, x_3} = 24$  is virtually indistinguishable from the analytical solution. Next, we visualize the resulting vertical boundary velocity $g_{x_2}$ and the corresponding analytical solution in \cref{tikz:boundvelres}. The velocities are linear within the element due to the $\mathrm{C}^1$ mesh continuity. Again, the results on the finest mesh with $n_{\mathrm{el}, x_3} = 24$ are virtually indistinguishable from the analytical solution.

\clearpage

\subsection{Stationary hydrofoil}
We simulate the flow past a stationary hydrofoil for angles of attack $\alpha$ ranging from 1$ ^{\circ}$ to 5$ ^{\circ}$. The simulations are performed on a NACA0012 foil section with Reynolds number $\mathbb{R}{\rm e} = U c/\nu$ = 1000 where $\nu$ is the kinematic viscosity. We study the resulting drag coefficient $C_d$ and lift coefficient $C_l$ defined as:
\begin{subequations}
\begin{alignat}{2}
    C_d =&~ \frac{2 F_d}{\rho c U^2 },\\
    C_l =&~ \frac{2 F_l}{\rho c U^2 },
\end{alignat}
\end{subequations}
where $F_d$ is the force component in the flow direction, $F_l$ the force component perpendicular to the flow direction, and $\rho$ denotes the density. 

\begin{figure}[h]
    \centering
    \begin{subfigure}[b]{0.49\textwidth}
        \centering
        \begin{tikzpicture}
    \begin{axis}[
        xmin = 0, xmax = 1,
        ymin = 0.121, ymax = 0.124,
        xtick distance = 0.25,
        ytick distance = 0.001,
        width = \textwidth,
        height = \textwidth,
        xlabel = {$h/h_0$ $[-]$},
        ylabel = {$C_d$ $[-]$},
        domain=0:3,
        legend pos = north west,
        ]
        \addplot +[only marks]table[x index=0,y index=1] {data/steady/2Dgridconv.dat};
        \addplot +[line width=1.0pt, mark=none, black, dotted]table[x index=0,y index=1] {data/steady/2Dgridconvrich.dat};
        \legend{2D, Rich. ext. $h/h_0 \xrightarrow{} 0$ }
        \end{axis}
    \end{tikzpicture}
    \caption{}
    \label{tikz:gridconvcd}
    \end{subfigure}
    \hfill
    \begin{subfigure}[b]{0.49\textwidth}
    \centering
\begin{tikzpicture}
\begin{axis}[
        xmin = 0, xmax = 1,
        ymin = 0.158, ymax = 0.161,
        xtick distance = 0.25,
        ytick distance = 0.001,
        width = \textwidth,
        height = \textwidth,
        xlabel = {$h/h_0$ $[-]$},
        ylabel = {$C_l$ $[-]$},
        domain=0:3,
        legend pos = north west,
        ]
        \addplot +[only marks]table[x index=0,y index=2] {data/steady/2Dgridconv.dat};
        \addplot +[line width=1.0pt, mark=none, black, dotted]table[x index=0,y index=2] {data/steady/2Dgridconvrich.dat};
        \legend{2D, Rich. ext. $h/h_0 \xrightarrow{} 0$ }
        \end{axis}
    \end{tikzpicture}
    \caption{}
    \label{tikz:gridconvcl}
    \end{subfigure}
    \caption{Results of steady 2D simulations for the drag coefficient $C_d$ in (a) and lift coefficient $C_l$ in (b) using four meshes and a Richardson extrapolation of the limit $h/h_0 \xrightarrow{} 0$ based on the three fines meshes. The order of convergence is 1.57 and 1.34 for drag and lift respectively. }
    \label{tikz:gridconv}
\end{figure}
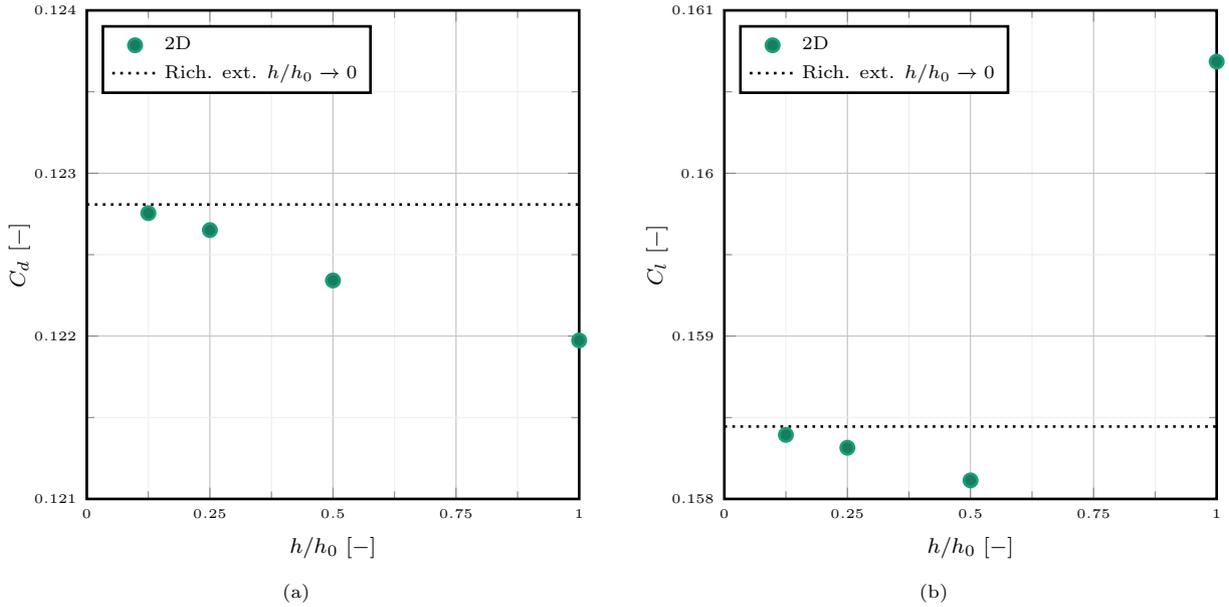

We first consider the two-dimensional setup. \cref{tikz:gridconv} shows the results of the spatial grid convergence study for $C_d$ and $C_l$ using 4 different meshes of varying resolution. In the coarsest mesh the domain is discretized using $30$ elements over the length of the hydrofoil, $15$ elements between the hydrofoil and the inflow boundary, and $45$ elements between the hydrofoil and the outflow boundary. We use a Richardson extrapolation to examine the limit $h/h_0 \to 0$ using the three finest meshes only, as the coarsest mesh is not in the asymptotic range. We find $1.57$ and $1.34$ for the order of convergence of the drag and lift, respectively. We choose the mesh with two refinements for our computations as this gives a balance between results and computational efforts. For this mesh the error is $0.13\%$ and $0.08 \%$ for $C_d$ and $C_l$ respectively considering the extrapolated result for $h/h_0 \to 0$. 

\clearpage

\begin{figure}[h]
    \centering
    \begin{tikzpicture}
        \begin{axis}[
        xmin = 0, xmax = 6,
        ymin = 0.0, ymax = 0.3,
        xtick distance = 1,
        ytick distance = 0.05,
        width =  0.5\textwidth,
        height = 0.5\textwidth,
        xlabel = {$\alpha$ $[^\circ]$},
        ylabel = {$C_l [-]$ },
        legend pos = north west,
        ]
        \addplot +[only marks, mark size = 3.75]table[x index=0,y index=1] {data/steady/steady2Daoa-delfi.dat};
        \addplot +[only marks, mark size = 2.25]table[x index=0,y index=1] {data/steady/steadySTaoa-delfi.dat};
        
        \addplot +[only marks, mark = +, black, line width=0.5pt]table[x index=0,y index=1] {data/steady/steady2Daoa-xfoil.dat};
        \addplot +[only marks, mark = x, black, line width=0.5pt]table[x index=0,y index=1]{data/steady/steady2Daoa-kurtulus.dat};
        \addplot +[only marks, mark = star, black, line width=0.5pt]table[x index=0,y index=1]{data/steady/steady2Daoa-liu.dat};
        \addplot +[only marks, mark = Mercedes star, black, line width=0.5pt]table[x index=0,y index=1]{data/steady/steady2Daoa-khalid.dat};
        \legend{ 2D, ST,  BEM \cite{Drela1989},  RANS \cite{Kurtulus2015}, ALE-CBS \cite{Liu2012}, RANS \cite{Khalid2012}}
        \end{axis}
        \end{tikzpicture}
    \caption{ Lift coefficient $C_l$ at $\mathbb{R}{\rm e} = 1000$ of a stationary NACA0012 hydrofoil for several angles of attack $\alpha$ determined using the proposed method supplemented with results from the literature. }
    \label{tikz:steadyaoa}
\end{figure}
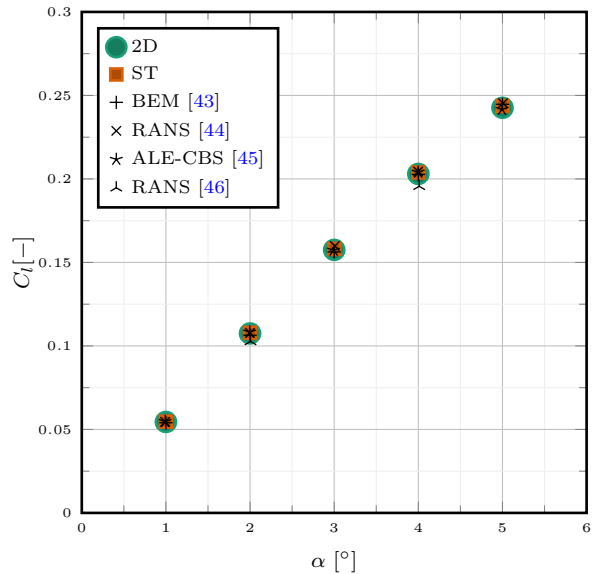

Next, we focus on the lift coefficient. \cref{tikz:steadyaoa} shows $C_l$ determined in space-time and two-dimensional simulations, supplemented with results from the literature. The computations are performed for $5$ different angles of attack. The similarity of the results of space-time and two-dimensional simulations demonstrates that the spatial convergence of two-dimensional simulations is indeed sufficient for space-time simulations. Moreover, the results are in good agreement with the results from the literature. We compare with (i) a Boundary Element Method (BEM) with viscous correction XFoil \cite{Drela1989}, (ii) the Reynolds Averaged Navier Stokes solver Ansys Fluent \cite{Kurtulus2015} (RANS), (iii) an Arbitrary-Lagrangian-Eulerian Characteristic Based Split Scheme solver \cite{Liu2012} (ALE-CBS), and (iv) Ansys Fluent computations \cite{Khalid2012}. The latter computations are only available for the angles of attack of 2$ ^{\circ}$ and 4$^{\circ}$. The numerical results obtained with this solver deviate more from the results that we have obtained. Lastly we note that we have verified the force signal of the space-time simulations to be constant in time. This demonstrates that our method correctly predicts steady flow.

\subsection{Heaving hydrofoil at a low reduced frequency}
In this test case we simulate a slowly heaving hydrofoil. The hydrofoil is oscillating at a low reduced frequency $k = {\pi c}/{(T U)}$. We note that the effect of the unsteady wake on the flow past the hydrofoil is very low \cite{Theodorsen1949} and added mass effects are negligible. As a consequence, the forces on the hydrofoil should match these from quasi-static simulations. We obtain the quasi-static results using two-dimensional simulations where we compensate the angle of attack $\alpha$ for inflow due to the heave motion. This provides the effective angle of attack:
\begin{equation}
    \alpha_{\mathrm{eff}} = \alpha - \arctan \left( \frac{ 2\pi h_a \cos\left(2 \pi T^{-1} t\right)}{T U}\right).
\end{equation}
The simulations are performed with $\mathbb{R}{\rm e} = 1000$, $k =$ 0.01, $h_a$ = 0.1 $m$ and $\alpha = 0^\circ$. We use the same spatial discretization as for the stationary cases. In the temporal direction we use $n_{\mathrm{el}, x_3} =  24$. We note that further refinement does not improve the numerical results.  

\begin{figure}
    \centering
        \begin{subfigure}[t]{0.49\textwidth}
        \vskip 0pt
        \centering
        \begin{tikzpicture}
\begin{axis}[
    ymode=log,
    xmin = 0, xmax = 55,
    ymax = 1000,
    xtick distance = 10,
    ytick distance = 100,
    width = \textwidth,
    height =\textwidth,
    xlabel = {$\theta$ $[s]$},
    ylabel = {$||\mathrm{\br}||_2$},
    legend pos = north east,
    ]
\addplot +[mark = none]table[x index=0, y index=2] {data/slowheave/slowheaveConvergence.dat};
\addplot +[mark = none]table[x index=0, y index=3] {data/slowheave/slowheaveConvergence.dat};
\addplot +[mark = none]table[x index=0, y index=4] {data/slowheave/slowheaveConvergence.dat};
\legend{mom.$_{x_1}$, mom.$_{x_2}$, mass}
\end{axis}
\end{tikzpicture}
    \caption{}
    \label{tikz:convergence}
    \end{subfigure}
    \hfill
    \begin{subfigure}[t]{0.49\textwidth}
        \vskip 0pt
        \centering
        \begin{tikzpicture}
\def\T{314.16}
\begin{axis}[
    xmin = 0, xmax = 1,
    ymin = -0.05, ymax = 0.15,
    xtick distance = 0.250,
    ytick distance = 0.05,
    width = \textwidth,
    height =\textwidth,
    xlabel = {$t/T$ $[-]$},
    ylabel = {$C_f$ $[-]$},
    domain=0:1,
    legend style={at={(0.03,0.5)},anchor=west},
    ]
\addplot +[mark = none]table[x expr=\thisrow{tx}/\T,y index=1] {data/slowheave/slowheave.dat};
\addplot +[mark = none]table[x expr=\thisrow{ty}/\T,y index=3] {data/slowheave/slowheave.dat};
\addplot [only marks, mark = x, black, line width=0.5pt]table[x expr=\thisrow{t}/\T,y index=1]
{data/slowheave/slowheave2D.dat};
\addplot +[only marks, mark = star, black, line width=0.5pt]table[x expr=\thisrow{t}/\T, y index=2]
{data/slowheave/slowheave2D.dat};
\legend{$C_d$ ST, $C_l$ ST,  $C_d$ QS, $C_l$ QS}
\end{axis}
\end{tikzpicture}
    \caption{}
    \label{tikz:slowheav}
    \end{subfigure}
    \caption{The convergence and results of the case of a slowly sinusoidal heaving NACA0012 hydrofoil with $k$ = 0.01 and $\mathbb{R}{\rm e}$ = 1000: a) $L^2$-norm of the residuals at the start of the first Newton iteration of space-time momentum and mass conservation over pseudo-time $\theta$; b) Force coefficients $C_f = C_d, C_l$ in space-time (ST) compared to quasi-static (QS) results.}
    \label{tikz:slowheavtot}
\end{figure}
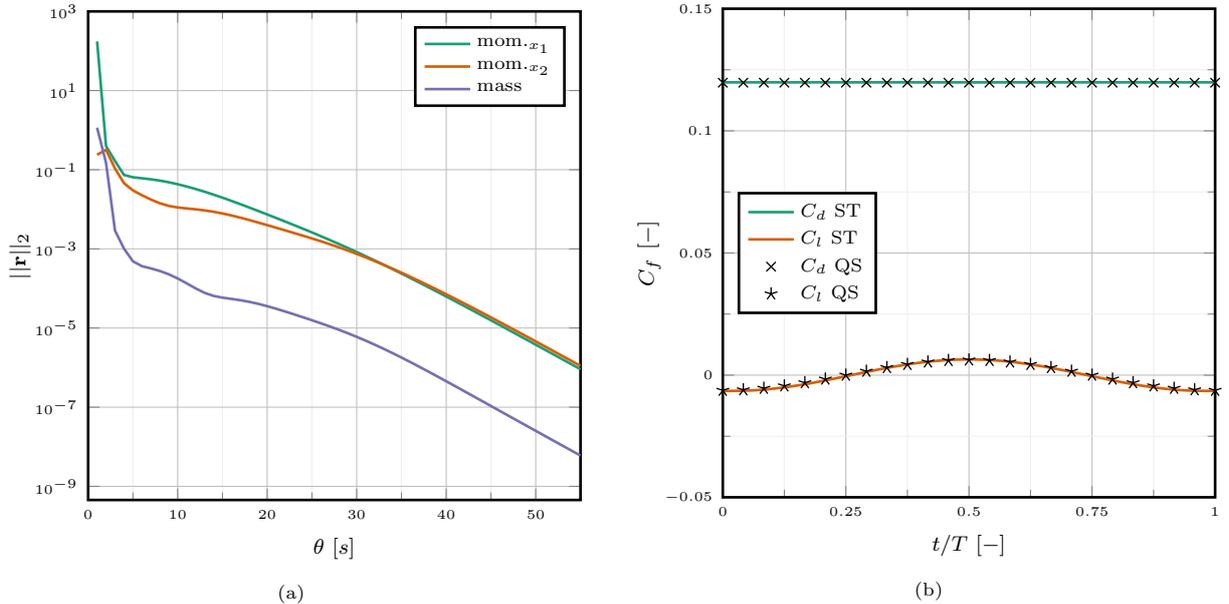

We visualize the convergence of the residuals in \cref{tikz:convergence}. We have verified that using stricter convergence criteria does not improve the solution quality. In \cref{tikz:slowheav} we show $C_l$ and $C_d$ for the space-time and quasi-static simulations. We observe that both $C_l$ and $C_d$ agree with the quasi-static results.

\subsection{Hydrofoil with large angle pitch motion}
In this last test case we focus on the prediction of the history effects in the wake. We simulate the flow past a sinusoidal pitching NACA0015 hydrofoil. The hydrofoil pitches around the \nicefrac{1}{3} chord with motion $\alpha(t) = \alpha_a \sin (2\pi t/T)$, where the amplitude is $\alpha_a = 23^\circ$, the Reynolds number is $\mathbb{R}{\rm e} = 1100$ and frequency is $k = 0.377$.

The same case is studied by \cite{Kinsey2008} using Ansys Fluent. Their simulation setup uses an impulsive start and at least $20$ large time steps to move the wake downstream of the hydrofoil. Their simulation is pursued with more than $2000$ time steps per period, and the simulation is considered periodic if the maximum variation in mean statistics between the last cycles is $0.1 \%$. In our setup we use the same spatial discretization as in our previous space-time simulations. To accurately capture the flow characteristics, we apply two extra refinements in the temporal direction. We note that further refinement does not yield improved solution quality. 

\begin{figure}
\begin{tikzpicture}
\def\T{8.33}
\begin{axis}[
    xmin = 0, xmax = 1,
    ymin = -1.5, ymax = 1.5,
    xtick distance = 0.2,
    ytick distance = 0.5,
    width = \textwidth,
    height = 0.5\textwidth,
    xlabel = {$t/T$ $[-]$},
    ylabel = {$C_l$ $[-]$},
    domain=0:\T,
    legend pos = north west,
    ]
\addplot +[mark = none]table[x expr=\thisrow{t}, y index=1] {data/kinsdumas/kinsey-dumasST-4-2.dat};
\legend{ ST, \cite{Kinsey2008} }
\addplot +[only marks, mark = x, black, line width=0.5pt]table[x expr=\thisrow{t}/\T,y index=1] {data/kinsdumas/kinsey-dumas.dat};
\end{axis}
\end{tikzpicture}
\caption{Lift coefficient $C_l$ over time of a pitching NACA0015 hydrofoil with $\alpha_a$ = 23$^\circ$, l and $\mathbb{R}{\rm e} = 1100$.}
\label{tikz:pitch}
\end{figure}
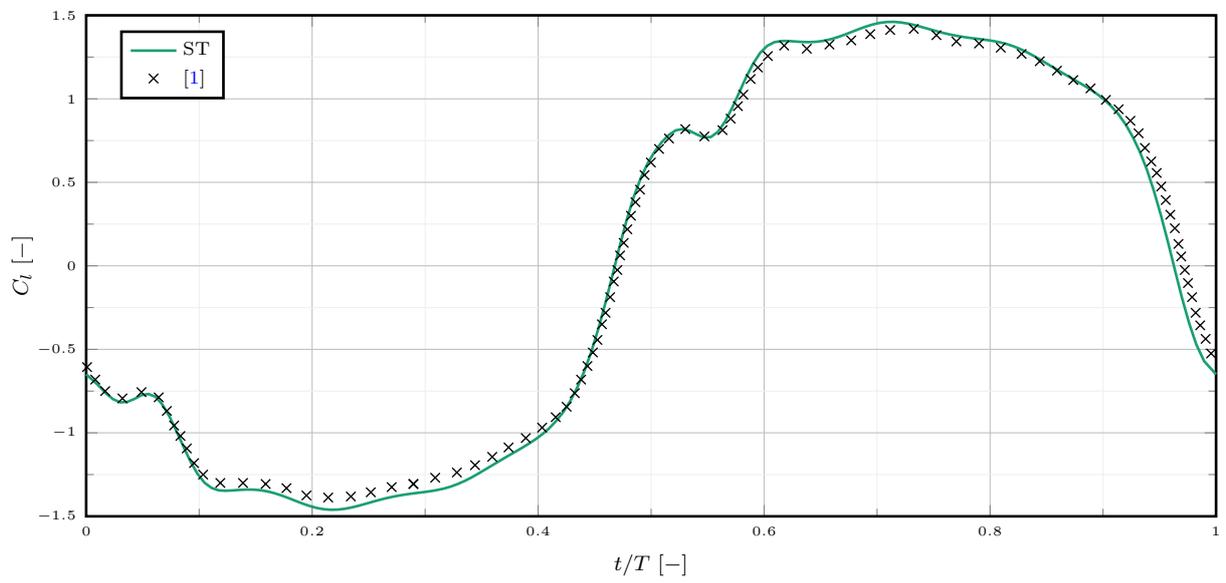

In \cref{tikz:pitch} we show a time signal of the lift coefficient $C_l$. In general we observe good agreement between our result and the result of \cite{Kinsey2008}. We see small differences in the regions $0.10 < t/T < 0.43$ and $0.58 < t/T < 0.84$. One important difference between our setup and the simulation in \cite{Kinsey2008} is that our solution is exactly periodic which is not the case in the reference computation. In \cref{fig:cfd} we show the velocity and pressure fields for $8$ moments in time. Note the periodic solution behavior. This is most apparent in the flow behind the hydrofoil when comparing the velocity field at $t/T = 7/8$ and $t/T = 0$. Furthermore, note that the flow is symmetric around the $x$-axis. To see this, compare for instance the velocity field at $t/T = 0$ with $t/T = 4/8$ and $t/T = 2/8$ with $t/T = 6/8$. 
Both figures illustrate that the effect of the history in the wake is correctly predicted.

\newpage
\newgeometry{margin = 0.5 cm}
\begin{landscape}
\centering

\begin{figure}
\thispagestyle{empty}
\centering
\begin{tikzpicture}
    % Insert main figure
    \node (image) at (0,0) {
        \includegraphics[width= 1.4\textheight]{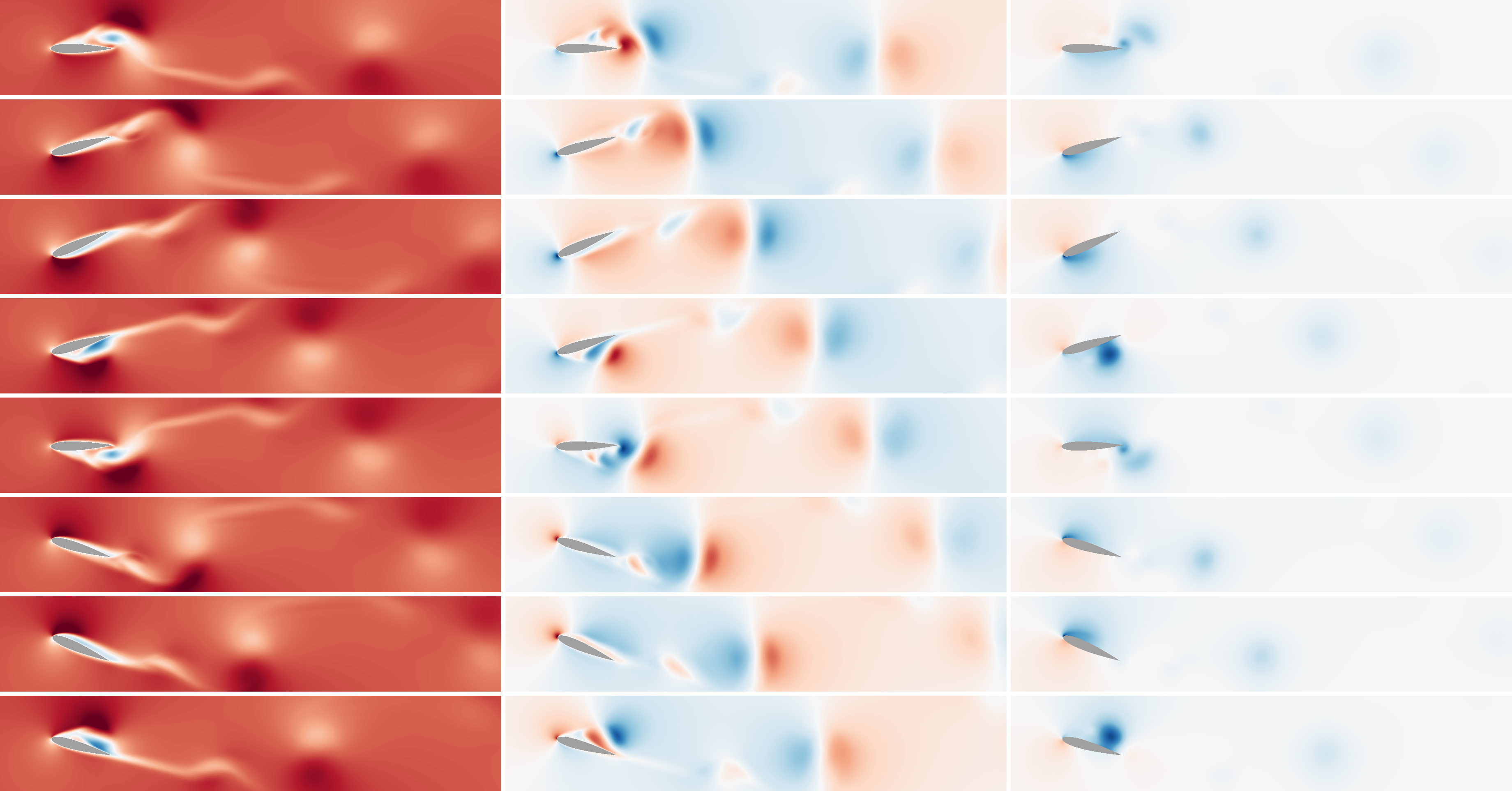}
    };
    
    % Time labels
    \node [fill=white,inner sep=1pt] at (-13.3,  7.15) {\scriptsize t/T = 0/8};
    \node [fill=white,inner sep=1pt] at (-13.3,  5.31) {\scriptsize t/T = 1/8};
    \node [fill=white,inner sep=1pt] at (-13.3,  3.46) {\scriptsize t/T = 2/8};
    \node [fill=white,inner sep=1pt] at (-13.3,  1.63) {\scriptsize t/T = 3/8};
    \node [fill=white,inner sep=1pt] at (-13.3, -0.22) {\scriptsize t/T = 4/8};
    \node [fill=white,inner sep=1pt] at (-13.3, -2.06) {\scriptsize t/T = 5/8};
    \node [fill=white,inner sep=1pt] at (-13.3, -3.91) {\scriptsize t/T = 6/8};
    \node [fill=white,inner sep=1pt] at (-13.3, -5.75) {\scriptsize t/T = 7/8};

    % Insert legends
        % Leg u_1
    \def\os{0}
    \node [fill=white, draw = black, inner sep=3pt] at (-13.65+\os+0.5, 7.75) {\Large  $u_1$ \normalsize [m/s]};
    \node (image) at (-6.1+\os-0.5,7.05) {\scalebox{-1}[1]{
        \includegraphics[scale=0.15,angle=90,origin=c]{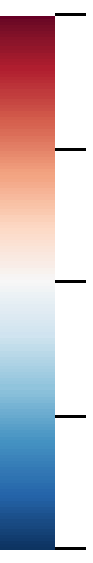}}
    };
    \node [fill=white,inner sep=1pt] at (-7.55+\os-0.5,8.0) {\scriptsize  -1.5};
    \node [fill=white,inner sep=1pt] at (-6.15+\os-0.5,8.0) {\scriptsize  0};
    \node [fill=white,inner sep=1pt] at (-4.75+\os-0.5,8.0) {\scriptsize  1.5};
    
        % Leg u_2
    \def\os{9.35}
    \node [fill=white, draw = black, inner sep=3pt] at (-13.65+\os+0.5, 7.75) {\Large  $u_2$ \normalsize [m/s]};
    \node (image) at (-6.1+\os-0.5,7.05) {\scalebox{-1}[1]{
        \includegraphics[scale=0.15,angle=90,origin=c]{fig/kins-dumas/legend.png}}
    };
    \node [fill=white,inner sep=1pt] at (-7.55+\os-0.5,8.0) {\scriptsize  -1.2};
    \node [fill=white,inner sep=1pt] at (-6.15+\os-0.5,8.0) {\scriptsize  0};
    \node [fill=white,inner sep=1pt] at (-4.75+\os-0.5,8.0) {\scriptsize  1.2};
    
        % Leg p
    \def\os{2*9.35}
    \node [fill=white, draw = black, inner sep=3pt] at (-13.65+\os+0.7, 7.75) {\Large  $p$ \normalsize [N m/kg]};
    \node (image) at (-6.1+\os-0.5,7.05) {\scalebox{-1}[1]{
        \includegraphics[scale=0.15,angle=90,origin=c]{fig/kins-dumas/legend.png}}
    };
    \node [fill=white,inner sep=1pt] at (-7.55+\os-0.5,8.0) {\scriptsize  -1.7};
    \node [fill=white,inner sep=1pt] at (-6.15+\os-0.5,8.0) {\scriptsize  0};
    \node [fill=white,inner sep=1pt] at (-4.75+\os-0.5,8.0) {\scriptsize  1.7};
    
\end{tikzpicture}
\caption{Velocity and pressure plots of a pitching hydrofoil with an angle of 23$^{\circ}$ and a period $T$ = 8.33 s for 8 moments in time. $\mathbb{R}{\rm e} = 1100$.}
\label{fig:cfd}
\end{figure}

\end{landscape}
\restoregeometry
\newpage

\section{Conclusions}
\label{sec:conclusions}
In this work we introduce a time-periodic continuous space-time formulation to simulate flow past periodically moving objects. The method employs isogeometric analysis to achieve higher-order smoothness in space and time. We discretize the formulation using residual-based variational turbulence modeling in which turbulent eddy viscosities are absent. Furthermore, we use weak boundary conditions to enhance the accuracy near the moving boundaries of the computational domain and pseudo-transient continuation to overcome some of the difficulties related to the saddle-point nature of the underlying problem. We show the conservation properties of the formulation and introduce a conservative traction evaluation. 
Numerical experiments on flow past stationary and moving hydrofoils demonstrate very good accuracy, even on coarse meshes. The computed drag and lift coefficients match results from literature and history effects in the wake are accurately captured. The proposed methodology circumvents lengthy transient simulations of periodic flow. A single period computation suffices, and the numerical solution is exactly periodic by design.

\section{Acknowledgements}
This publication is part of the project Lift Control for Hydrofoil Craft (with project number TWM. BL.019.009) of the research programme Top Sector Water \& Maritime: The Blue Route which is (partly) financed by the Dutch Research Council (NWO). This support is gratefully acknowledged. MtE acknowledges support from the German Research Foundation (Deutsche Forschungsgemeinschaft DFG) via the Walter Benjamin project EI 1210/1-1.

\appendix
%\section*{References}

\bibliography{mybibfile}

\end{document}